\documentclass[12pt,twoside,a4paper]{article}
\usepackage{enumerate}
\usepackage{t1enc} 
\usepackage{times}
\usepackage[dvips]{epsfig}
\usepackage{amssymb}
\usepackage{amsmath}
\usepackage{longtable}
\usepackage{makeidx}
\usepackage{euscript}
\usepackage{enumerate}
\usepackage{hyperref}
\DeclareMathAlphabet\scr{U}{scr}{m}{n}
\SetMathAlphabet\scr{bold}{U}{scr}{b}{n}
  \DeclareFontFamily{U}{scr}{\skewchar\font'177}%
  \DeclareFontShape{U}{scr}{m}{n}{<-6>rsfs5<6-8>rsfs7<8->rsfs10}{}%
  \DeclareFontShape{U}{scr}{b}{n}{<-6>rsfs5<6-8>rsfs7<8->rsfs10}{}%

\newcommand{\rr}{\mathbb R}  
\newcommand{\ev}{\mathbb E}

\newcommand{\nn}{\mathbb N}

\newcommand{\pp}{\mathbb P}

\newtheorem{satz}{Theorem}[section]
\newtheorem{prop}[satz]{Proposition}
\newtheorem{theorem}[satz]{Theorem}

\newtheorem{lemma}[satz]{Lemma}

\newtheorem{assumption}[satz]{Assumption}

\newtheorem{@definition}[satz]{Definition}

\newtheorem{@aufgabe}{Aufgabe}

\newtheorem{@bsp}[satz]{Beispiel}

\newenvironment{mat}{\left (\begin{matrix} } {\end{matrix}\right ) }

\newcommand{\bmat}{\begin{mat}}
\newcommand{\emat}{\end{mat}}
\newcommand{\be}{\begin{enumerate}}
\newcommand{\ee}{\end{enumerate}}
\newcommand{\beq}{\begin{equation}}
\newcommand{\eeq}{\end{equation}}
\newcommand{\bea}{\begin{eqnarray}}
\newcommand{\eea}{\end{eqnarray}}
\newcommand{\beaa}{\begin{eqnarray*}}
\newcommand{\eeaa}{\end{eqnarray*}}

\newcommand{\ep}{\hfill $\Box$}

\renewcommand{\epsilon}{\varepsilon}
\renewcommand{\phi}{\varphi}
\renewcommand{\rho}{\varrho}

\title{{\small\bf LARGE VOLATILITY-STABILIZED MARKETS}\footnote{Research supported in part by NSF grant DMS-0806211.\newline
{\it AMS 2000 subject classifications:} Primary 60J60; secondary 60J70, 91G80.\newline
{\it Keywords and phrases:} interacting diffusion processes, hydrodynamic limit, volatility-stabilized models, Bessel processes, degenerate parabolic partial differential equations.}}
\author{{\small BY M. SHKOLNIKOV}\\\quad\\{\small\it Stanford University}}
\date{}

\begin{document}
\maketitle

\begin{abstract}
We investigate the behavior of systems of interacting diffusion processes, known as volatility-stabilized market models in the mathematical finance literature, when the number of diffusions tends to infinity. We show that, after an appropriate rescaling of the time parameter, the empirical measure of the system converges to the solution of a degenerate parabolic partial differential equation. A stochastic representation of the latter in terms of one-dimensional distributions of a time-changed squared Bessel process allows us to give an explicit description of the limit. 
\end{abstract}

\section{Introduction}

Recently, Fernholz and Karatzas \cite{fk} have introduced two types of systems of interacting diffusion processes, the volatility-stabilized market models and the rank-based market models, in the context of stochastic portfolio theory. Both of them serve as models for the evolution of capitalizations in large financial markets and incorporate the fact that stocks of firms with smaller market capitalization tend to have higher rates of returns and be more volatile. In a previous paper \cite{sh} the author gave a description of the joint dynamics of the market capitalizations in rank-based models, when the number of firms tends to infinity (see also \cite{jo} for related results). Here, the corresponding limit is investigated in the context of volatility-stabilized models.  

\bigskip

The dynamics of the capitalizations in volatility-stabilized models is given by the unique weak solution to the following system of stochastic differential equations:
\begin{eqnarray}\label{sde1}
dX_i(t)=\frac{\eta}{2}S(t)dt+\sqrt{X_i(t)S(t)}\;dW_i(t),\quad 1\leq i\leq N,
\end{eqnarray}
which is endowed with an initial distribution of the vector $(X_1(0),\dots,X_N(0))$ on $[0,\infty)^N$. Hereby, $\eta$ is a real number greater than $1$, $S(t)=X_1(t)+\dots+X_N(t)$ and $W_1,\dots,W_N$ is a collection of $N$ independent standard Brownian motions. We refer to section 12 of \cite{fk} for a construction of a weak solution to \eqref{sde1} and an explanation why it is unique.  

\bigskip

We will analyze the limit of the path of empirical measures $\frac{1}{N}\sum_{i=1}^N \delta_{X_i(t)}$ corresponding to \eqref{sde1}, after a suitable rescaling of the time parameter, when $N$ tends to infinity. The slowdown of the time by a factor of $N$ is needed to observe a non-degenerate limiting behavior. Heuristically, this can be inferred from the appearance of the process $S(t)=X_1(t)+\dots+X_N(t)$, an order $N$ object, in the drift and diffusion coefficients of the processes $X_1,\dots,X_N$. 

\bigskip 
     
We show that the limit of the sequence of laws of $\frac{1}{N}\sum_{i=1}^N \delta_{X_i(t/N)}$, $N\in\nn$ exists and that under the limiting measure the degenerate linear parabolic equation \eqref{pde} below must be satisfied in the weak sense with probability $1$. Using a stochastic representation of the solution to \eqref{pde}, we can determine the latter explicitly. Hence, our results allow to approximate the evolution of the capitalizations in a large volatility-stabilized market by the solution of the limiting equation \eqref{pde}. Moreover, in the context of stochastic portfolion theory (see e.g. \cite{fe}, \cite{fk}) one is interested in the behavior of the rank statistics of the vector $(X_1(t),\dots,X_N(t))$ of capitalizations. Since these are given by the $\frac{1}{N},\frac{2}{N},\dots,\frac{N}{N}$-quantiles of the empirical measure $\frac{1}{N}\sum_{i=1}^N \delta_{X_i(t)}$, our results can be also used to approximate the evolution of any finite number of ranked capitalizations by the evolution of the corresponding quantiles of the solution to the partial differential equation \eqref{pde}. In addition, the stochastic representation mentioned above shows that the solution to the equation \eqref{pde} is given by the one-dimensional distributions of a time-changed squared Bessel process and, thus, establishes a new connection between volatility-stabilized market models and squared Bessel processes (see \cite{fk} for further connections). The latter were analyzed in much detail in the works \cite{py1}, \cite{py2} and \cite{ry} among others. 

\bigskip

Independently from the field of stochastic portfolio theory, systems of interacting diffusion processes play a major role in statistical physics. In particular, systems of diffusions interacting through their empirical measure (mean field) were studied in the literature by many authors, see e.g. \cite{ga}, \cite{mc}, \cite{da}, \cite{fu}, \cite{le}, \cite{nt1}, \cite{nt2}, \cite{oe}. We remark that the system \eqref{sde1} can be cast into the framework of \cite{ga}, since the drift and the diffusion coefficients in the $i$-th equation of the system \eqref{sde1} can be expressed as functions of the empirical measure of the particle system and the position of the $i$-th particle. However, the generator of the particle system is not uniformly elliptic on $[0,\infty)^N$ and the same is true on $[0,\infty)$ for the elliptic differential operator on the right-hand side of the equation \eqref{pde}. For this reason, the results of \cite{ga} do not carry over directly to our setting. Nonetheless, we adapt some of the techniques developed there to our case.

\bigskip   

The time-varying mass partition 
\begin{eqnarray}
\alpha_i(t)=\frac{X_i(t)}{X_1(t)+\dots+X_N(t)},\quad 1\leq i\leq N
\end{eqnarray}
is referred to as the collection of market weights in the mathematical finance literature and describes the capitalizations of the firms as fractions of the total capitalization of the market. The collection of market weights has the remarkable property of being extremely stable over time for all major financial markets (see \cite{fe} for plots of the market weights in multiple real-world markets). This was explained to a large extent in the context of rank-based market models in \cite{cp} and \cite{ps}. 

\bigskip

The model \eqref{sde1} incorporates the empirically observed fact that the capitalizations of firms with a small market weight tend to have a higher rate of growth and to fluctuate more wildly. Indeed, this becomes apparent from the dynamics of the logarithmic capitalizations corresponding to \eqref{sde1}:
\begin{eqnarray}
d(\log X_i(t))=\frac{\eta-1}{2\alpha_i(t)}\;dt+\frac{1}{\sqrt{\alpha_i(t)}}\;dW_i(t),\quad 1\leq i\leq N
\end{eqnarray}
and the assumption $\eta>1$. A detailed analysis of the evolution of the market weights under the model \eqref{sde1} and the corresponding invariant measure can be found in \cite{pa}.

\bigskip

We assume the following condition on the initial values $X_1(0),\dots,X_N(0)$ of the capitalizations. 

{\assumption\label{ass} The laws of $\frac{1}{N}\sum_{i=1}^N \delta_{X_i(0)}$, $N\in\nn$ on $M_1([0,\infty))$, the space of probability measures on $[0,\infty)$ endowed with the topology of weak convergence of measures, converge weakly to $\delta_\lambda$ for some $\lambda\in M_1([0,\infty))$ with a finite first moment, the quantities $\ev[S(0)]$ and $\ev[S(0)^2]$ are finite for all $N\in\nn$, and it holds
\begin{eqnarray}\label{m_lambda}
\lim_{N\rightarrow\infty}\frac{\ev[S(0)]}{N}=m_\lambda,\quad\lim_{N\rightarrow\infty}\frac{\ev[S(0)^2]}{N^2}=m_\lambda^2,
\end{eqnarray} 
where $m_\lambda=\int_{[0,\infty)} x\;\lambda(dx)$. In addition, we make the non-degeneracy assumption $m_\lambda>0$.}

\bigskip

We remark at this point that Assumption \ref{ass} is, in particular, satisfied if the random variables $X_1(0),\dots,X_N(0)$ are i.i.d. and distributed according to a measure $\lambda$ with finite two first moments and $m_\lambda>0$. This is a consequence of Varadarajan's Theorem in the form of Theorem 11.4.1 in \cite{du2}. 

\bigskip 

In order to formulate our main results we introduce the following set of notations. We write $M_1(\rr)$ and $M_1([0,\infty))$ for the spaces of probability measures on the real line and the non-negative half-line, respectively. We metrize both spaces in a way compatible with the topology of weak convergence of measures. Moreover, for a positive real number $T$, we let $C([0,T],M_1(\rr))$ and $C([0,T],M_1([0,\infty)))$ be the spaces of continuous functions from $[0,T]$ to $M_1(\rr)$ and from $[0,T]$ to $M_1([0,\infty))$, respectively, endowed with the topology of uniform convergence. In addition, we introduce the time-changed capitalization processes $Y_i(t)=X_i(t/N)$, $t\geq0$, $1\leq i\leq N$ and the corresponding path of empirical measures $\rho^N(t)=\frac{1}{N}\sum_{i=1}^N \delta_{Y_i(t)}$, $t\in[0,T]$ on an arbitrary finite time interval $[0,T]$, which is considered to be fixed from now on. Finally, we let $Q^N_T$ be the distribution of the random variable $\rho^N(t)$, $t\in[0,T]$ on $C([0,T],M_1(\rr))$. 

\bigskip

The main results of this paper are summarized in the following theorem.

{\theorem\label{main_result} Under Assumption \ref{ass} the following statements are true.
\begin{enumerate}[(a)]
\item The sequence $Q^N_T$, $N\in\nn$ converges weakly to a limit $Q^\infty_T$. Moreover, $Q^\infty_T$ is a Dirac delta measure, whose unique atom $\rho$ is given by the unique distributional solution of the Cauchy problem 
\begin{eqnarray}
&&\frac{\partial\rho}{\partial t}=-\frac{\eta}{2}e^{\frac{\eta t}{2}}m_\lambda\frac{\partial\rho}{\partial x}
+\frac{1}{2}e^{\frac{\eta t}{2}}m_\lambda \frac{\partial^2(x\rho)}{\partial x^2},\label{pde}\\
&&\rho(0)=\lambda \label{pdeic}
\end{eqnarray}
in $C([0,T],M_1([0,\infty)))$, where $m_\lambda=\int_{[0,\infty)} x\;\lambda(dx)$.
\item Let $Z(t)$, $t\geq0$ be a squared Bessel process satisfying the stochastic initial value problem
\begin{eqnarray}
&&dZ(t)=\frac{\eta}{2}\;dt+\sqrt{Z(t)}\;d\beta(t),\;t\geq0,\\
&&{\mathcal L}(Z(0))=\lambda,
\end{eqnarray}
where $\beta$ is a standard Brownian motion and ${\mathcal L}(Z(0))$ denotes the law of $Z(0)$. Then the unique distributional solution to the Cauchy problem \eqref{pde}, \eqref{pdeic} in $C([0,T],M_1([0,\infty)))$ is given by the one-dimensional distributions of the time-changed process $Z\Big(\int_0^t e^{\eta s/2}m_\lambda\;ds\Big)$, $t\in[0,T]$.
\item Let $\rho(t)$, $t\in[0,T]$ be the only atom of the measure $Q^\infty_T$. Then for every $t\in(0,T]$ the measure $\rho(t)$ is absolutely continuous with respect to the Lebesgue measure $Leb$ on $[0,\infty)$ and the corresponing density is given by  
\begin{eqnarray*}
\frac{d\rho(t)}{dLeb}(y)=\int_{[0,\infty)} \frac{2}{J(t)}\Big(\frac{y}{x}\Big)^{(\eta-1)/2}\exp\Big(-\frac{2(x+y)}{J(t)}\Big)I_{\eta-1}\Big(\frac{4\sqrt{xy}}{J(t)}\Big)\lambda(dx),
\end{eqnarray*}
where $J(t)=\int_0^t e^{\eta s/2}m_\lambda\;ds$ and $I_{\eta-1}$ is the Bessel function of the first kind of index $\eta-1$.
\end{enumerate}}

\bigskip

\noindent{\bf Remarks.} 
\begin{enumerate}[(1)]
\item By a distributional solution of the problem \eqref{pde}, \eqref{pdeic} in $C([0,T],M_1([0,\infty)))$ we mean an element $\rho$ of $C([0,T],M_1([0,\infty)))$ which satisfies the system \eqref{inteq2}, \eqref{inteq2ic}, where ${\mathcal S}(\rr)$ is the space of Schwartz functions on $\rr$ and $(\gamma,f)$ denotes $\int_\rr f\;d\gamma$ for any $f\in {\mathcal S}(\rr)$, $\gamma\in M_1(\rr)$. 
\item In \cite{sh} the limiting dynamics was derived (under some assumptions) for a different class of interacting diffusion processes, which go by the name of rank-based models in the context of stochastic portfolio theory. There, the limiting equation for the cumulative distribution function of the empirical measure of the logarithmic capitalizations was given by the porous medium equation, that is, a nonlinear non-degenerate parabolic partial differential equation. It was also shown there that its weak solution $w$ can be represented by the one-dimensional distributions of the process with the dynamics
\begin{eqnarray}
dX(t)=\mu(w(t,X(t)))\;dt+\sigma(w(t,X(t)))\;d\beta(t),
\end{eqnarray} 
where $\mu$ and $\sigma$ are functions depending on the parameters of the model. In contrast, the partial differential equation \eqref{pde} is a linear degenerate parabolic differential equation, which admits the stochastic representation of Theorem \ref{main_result} (b). Hence, although the rank-based and the volatility-stabilized models share multiple common features (such as the monotone dependence of the drift and diffusion coefficients of the logarithmic capitalizations on the market weights), their limiting behavior differs significantly. Indeed, it was shown in \cite{ar} that the weak solution of the porous medium equation may fail to be differentiable in the spatial variable for all $t\geq0$ in contrast to the findings in Theorem \ref{main_result} (c) for the equation \eqref{pde}. In addition, a big difference from the applicational point of view is the explicitness of the solution to the equation \eqref{pde}, whereas (in general) no explicit formula for the weak solution of the porous medium equation is known. It is also remarkable that with the usual parametrizations of the two models as in \cite{fk}, the time in the volatility-stabilized models \eqref{sde1} has to be slowed down by a factor of $N$ to observe non-degenerate limiting behavior, whereas this is not the case in the rank-based market models.  
\item Let $d$ be a metric on $M_1(\rr)$ which metrizes the topology of weak convergence of probability measures (such as the Levy metric or any other metric in section 11.3 of \cite{du2}), and let $d_{[0,T]}$ be the metric on $C([0,T],M_1(\rr))$ given by
\begin{eqnarray}
d_{[0,T]}(\xi_1,\xi_2)=\sup_{t\in[0,T]} d(\xi_1(t),\xi_2(t)).
\end{eqnarray}
Then $d_{[0,T]}$ makes the space $C([0,T],M_1(\rr))$ into a separable metric space (see e.g. Theorem 2.4.3 in \cite{sr}). Moreover, combining part (a) of Theorem \ref{main_result} with problem 6 in chapter 9.3 of \cite{du2}, we conclude that the sequence $\rho^N$, $N\in\nn$ converges to the path of measures $\rho$ of Theorem \ref{main_result} (c) in probability on $(C([0,T],M_1(\rr)),d_{[0,T]})$, that is,
\begin{eqnarray}
\forall\epsilon>0:\;\lim_{N\rightarrow\infty}\pp(\sup_{t\in[0,T]} d(\rho^N(t),\rho(t))>\epsilon)=0.
\end{eqnarray}
This gives an alternative way of stating part (a) of Theorem \ref{main_result}.   
\item The transition densities of the squared Bessel process $Z$ in Theorem \ref{main_result} (b) are known (see e.g. Corollary 1.4 in chapter XI of \cite{ry}), so that part (c) of Theorem \ref{main_result} is a direct consequence of parts (a) and (b) of Theorem \ref{main_result}.
\end{enumerate}

The rest of the paper is organized as follows. Part (a) of Theorem \ref{main_result} is proven in sections 2.1 and 2.2. Its proof is divided into three parts. Firstly, in Proposion \ref{tightness} it is shown that the sequence $Q^N_T$, $N\in\nn$ is tight. Its proof relies on the characterization of compact subsets of $C([0,T],M_1(\rr))$ obtained in \cite{ga} and a characterization of compact subsets of $C([0,T],\rr)$, the space of continuous real-valued functions on $[0,T]$ endowed with the topology of uniform convergence, given in \cite{sv}. Secondly, in Proposition \ref{convergence} we prove that under every limit point of the sequence $Q^N_T$, $N\in\nn$ the Cauchy problem \eqref{pde}, \eqref{pdeic} is satisfied in the distributional sense with probability $1$. Here, we use arguments from the theory of convergence of semimartingales in the spirit of \cite{js}. The main challenge in these two parts is to deal with the unboundedness of the drift and diffusion coefficients in the dynamics of the processes $Y_1,\dots,Y_N$. Thirdly, in Proposition \ref{uniqueness} we demonstrate that the problem \eqref{pde}, \eqref{pdeic} has a unique distributional solution in the space $C([0,T],M_1([0,\infty)))$. This is achieved by transforming the uniqueness problem into an existence problem and by applying existence results on boundary value problems for non-degenerate linear parabolic equations on bounded subsets of $[0,T]\times\rr$, on which the differential operator of equation \eqref{pde} is uniformly parabolic. After that, we give the proof of part (b) of Theorem \ref{main_result} in section 2.3 using methods of stochastic calculus. 

\section{Law of large numbers}
\setcounter{equation}{0}
\subsection{Tightness}

In this subsection we will combine the characterization of relatively compact sets in the space $C([0,T],M_1(\rr))$ of \cite{ga} with a characterization of relatively compact subsets of $C([0,T],\rr)$ in \cite{sv} to prove the tightness of the sequence $Q^N_T$, $N\in\nn$.

{\prop\label{tightness} The sequence $Q^N_T$, $N\in\nn$ is tight on $C([0,T],M_1(\rr))$.\\\quad\\}
{\it Proof.} 1) Let $C_c(\rr)$ be the space of compactly supported continuous functions on $\rr$ endowed with the topology of uniform convergence. We fix an arbitrary $\epsilon>0$ and a countable dense subset $\{f_1,f_2,\dots\}$ of $C_c(\rr)$ such that each $f_r$ is twice continuously differentiable. Moreover, for every $\gamma\in M_1(\rr)$ and every function $f$ on $\rr$, which is integrable with respect to $\gamma$, we write $(\gamma,f)$ for $\int_\rr f\;d\gamma$. 

From the proof of Lemma 1.3 in \cite{ga} we see that it is enough to find a compact set $K_0$ in $M_1(\rr)$ and compact sets $K_1,K_2,\dots$ in $C([0,T],\rr)$ such that for all $N\in\nn$:
\begin{eqnarray}
&&Q^N_T(\{\xi\in C([0,T],M_1(\rr))|\forall t\in[0,T]:\;\xi(t)\in K_0\})\geq 1-\epsilon,\\
&&Q^N_T(\{\xi\in C([0,T],M_1(\rr))|(\xi(.),f_r)\in K_r\})\geq 1-\epsilon\cdot 2^{-r},\; r\geq1. 
\end{eqnarray}  

To define $K_0$ we introduce the function $\phi(x)=|x|$ and use the non-negativity of the processes $Y_1,\dots,Y_N$ together with the dynamics of the processes $X_1,\dots,X_N$ to conclude
\begin{eqnarray*}
d(\rho^N(t),\phi)=\frac{\eta S^Y(t)}{2N} dt + \frac{1}{N^{3/2}}\sum_{i=1}^N \sqrt{Y_i(t)S^Y(t)} dB_i(t),
\end{eqnarray*}
where $B_i(t)=N^{1/2} W_i(t/N)$, $1\leq i\leq N$, $S^Y(t)=Y_1(t)+\dots+Y_N(t)$. From the representation of $X_1,\dots,X_N$ as time-changed squared Bessel processes (see equations (12.7)-(12.9) in \cite{fk}) and remark (ii) after Corollary 1.4 in chapter XI of \cite{ry} (note that their dimension parameter $\delta$ corresponds to our $2\eta$) it follows that the return time to $0$ of the processes $X_1,\dots,X_N$ is infinity with probability $1$. Hence, the process $B(t)=\sum_{i=1}^N \int_0^t \frac{\sqrt{Y_i(s)}}{\sqrt{S^Y(s)}} dB_i(s)$, $t\geq0$ is well-defined and a standard Brownian motion by Levy's Theorem. As a consequence we have
\begin{eqnarray}
d(\rho^N(t),\phi)=\frac{\eta}{2}(\rho^N(t),\phi) dt + \frac{(\rho^N(t),\phi)}{\sqrt{N}} dB(t).
\end{eqnarray}
The latter equation is a Black-Scholes stochastic differential equation and it is well-known that its unique strong solution is given by
\begin{eqnarray}
\;\;\;\;\;\;\;\;\;(\rho^N(t),\phi)=(\rho^N(0),\phi)\exp\Big((\eta/2-(2N)^{-1})t+N^{-1/2} B(t)\Big).
\end{eqnarray}  
Thus, for every $C>0$ and all $N\in\nn$ we have
\begin{eqnarray*}
&&\pp\Big((\sup_{0\leq t\leq T} (\rho^N(t),\phi))>C\Big)\\
&\leq&\pp\Big((\rho^N(0),\phi)\exp(N^{-1/2}\sup_{0\leq t\leq T} B(t))>C\exp(-(\eta/2-(2N)^{-1})T)\Big).
\end{eqnarray*}
A routine computation involving Chebyshev's inequality and \eqref{m_lambda} shows that the sequence of random variables $(\rho^N(0),\phi)\exp(N^{-1/2}\sup_{0\leq t\leq T} B(t))$, $N\in\nn$ converges in probability to the constant $m_\lambda$. Hence, by choosing $C$ large enough, one can make the latter upper bound smaller than $\epsilon$ for all $N\in\nn$. Thus, we can let $K_0$ be the closure of the set 
\begin{eqnarray*}
\{\gamma\in M_1(\rr)|(\gamma,\phi)\leq C\} 
\end{eqnarray*}
in $M_1(\rr)$, which is compact by Prokhorov's Theorem.

\bigskip

\noindent 2) To prove the existence of the sets $K_1,K_2,\dots$ with the desired properties it suffices to show that for any fixed $r\in\nn$ the sequence of probability measures $Q^{N,f_r}_T$, $N\in\nn$ on $C([0,T],\rr)$ induced by $Q^N_T$, $N\in\nn$ through the mapping $\xi\mapsto(\xi(.),f_r)$ is tight. To this end, we fix an $r\in\nn$ and aim to deduce the tightness of the sequence $Q^{N,f_r}_T$, $N\in\nn$ from Theorem 1.3.2 of \cite{sv}. To do this, we need to show 
\begin{eqnarray*}
\lim_{\theta\uparrow\infty}\inf_{N\in\nn} Q^{N,f_r}_T(|y(0)|\leq\theta)=1
\end{eqnarray*} 
and
\begin{eqnarray*}
\forall\Delta>0:\;
\lim_{\zeta\downarrow0}\limsup_{N\rightarrow\infty} Q^{N,f_r}_T\Big(\sup_{(s,t)\in A_\zeta}|y(t)-y(s)|>\Delta\Big)=0,
\end{eqnarray*}
where $A_\zeta=\{(s,t)|0\leq s\leq t\leq T,t-s\leq\zeta\}$. The first assertion follows immediately by considering $\theta>\sup_{x\in\rr}|f_r(x)|$. 

To prove the second assertion, we first rewrite it in terms of $X_1,\dots,X_N$:
\begin{eqnarray*}
\forall\Delta>0:\;\lim_{\zeta\downarrow0}\limsup_{N\rightarrow\infty} 
\pp\Big(\sup_{(s,t)\in A_{\zeta,N}}\frac{1}{N}\Big|\sum_{i=1}^N (f_r(X_i(t))-f_r(X_i(s)))\Big|>\Delta\Big)=0,
\end{eqnarray*}
where $A_{\zeta,N}=\{(s,t)|0\leq s\leq t\leq T/N,t-s\leq\zeta/N\}$. Next, we apply Ito's formula to the process $\frac{1}{N}\sum_{i=1}^N f_r(X_i(t))$, $t\geq0$ and conclude that it holds 
\begin{eqnarray}
\frac{1}{N}\sum_{i=1}^N f_r(X_i(t))=D(t)+M(t),\;t\geq0,
\end{eqnarray} 
where 
\begin{eqnarray}
\;\;\;\;\;\;\;\;\;D(t)&=&\frac{1}{N}\sum_{i=1}^N \int_0^t \frac{\eta}{2}f_r'(X_i(u))S(u)+\frac{1}{2}f_r''(X_i(u))X_i(u)S(u)\;du,\\
\;\;\;\;\;\;\;\;\;M(t)&=&\frac{1}{N}\sum_{i=1}^N \int_0^t f_r'(X_i(u))\sqrt{X_i(u)S(u)}\;dW_i(u).
\end{eqnarray}
Thus, for every fixed $\Delta>0$ the latter limit can be bounded from above by
\begin{eqnarray*}
&&\lim_{\zeta\downarrow0}\limsup_{N\rightarrow\infty}\pp\Big(\sup_{0\leq s\leq t\leq T/N,t-s\leq\zeta/N}|D(t)-D(s)|>\Delta/2\Big)\\
&+&\lim_{\zeta\downarrow0}\limsup_{N\rightarrow\infty}\pp\Big(\sup_{0\leq s\leq t\leq T/N,t-s\leq\zeta/N}|M(t)-M(s)|>\Delta/2\Big),
\end{eqnarray*}
which we will call expression (*). We will show that the first summand in (*) is zero in step 3 and that the second summand is equal to zero in step 4.

\bigskip

\noindent 3) To bound the first summand from above, we set 
\begin{eqnarray}
R=\max(\eta/2\sup_{x\in\rr}|f'_r(x)|,1/2\sup_{x\in\rr}|f''_r(x)|) 
\end{eqnarray}
and use the definition of the process $D$ to make the estimates
\begin{eqnarray*}
&&\pp\Big(\sup_{0\leq s\leq t\leq T/N,t-s\leq\zeta/N}|D(t)-D(s)|>\Delta/2\Big)\\
&\leq&\pp\Big(\sup_{0\leq s\leq t\leq T/N,t-s\leq\zeta/N}\int_s^t S(u)+S(u)^2/N\;du\geq\frac{\Delta}{2R}\Big)\\
&\leq&\pp\Big(\sup_{0\leq t\leq T/N}(S(t)+S(t)^2/N)\geq\frac{\Delta N}{2R\zeta}\Big)\\
&=&\pp\Big(\sup_{0\leq t\leq T/N} S(t)\geq\frac{-1+\sqrt{1+2\Delta/(R\zeta)}}{2}N\Big).
\end{eqnarray*}  

To estimate the latter upper bound further, we use the dynamics of the processes $X_1,\dots,X_N$ to find
\begin{eqnarray*}
dS(t)=S(t)\frac{\eta N}{2} dt + \sum_{i=1}^N \sqrt{X_i(t)S(t)} dW_i(t)
=S(t)\frac{\eta N}{2} dt + S(t) d\widetilde{B}(t),
\end{eqnarray*}  
where $\widetilde{B}(t)=\sum_{i=1}^N \int_0^t \frac{\sqrt{X_i(s)}}{\sqrt{S(s)}}dW_i(s)$, $t\geq0$ is a standard Brownian motion due to the same argument as for the process $B$ in step 1 of this proof. Thus, $S$ satisfies the Black-Scholes stochastic differential equation and is given explicitly by
\begin{eqnarray}\label{Sexplicitly}
S(t)=S(0)\exp\Big((\eta N/2-1/2)t+\widetilde{B}(t)\Big),\;t\geq0.
\end{eqnarray}

Hence, the latter upper bound is not greater than
\begin{eqnarray*}
\pp\Big(\frac{S(0)}{N}\exp\Big((\eta N/2-1/2)T/N+\sup_{0\leq t\leq T/N}\widetilde{B}(t)\Big)\geq \frac{-1+\sqrt{1+2\Delta/(R\zeta)}}{2}\Big).
\end{eqnarray*}
From \eqref{m_lambda} and Chebyshev's inequality it follows that the sequence of random variables 
\begin{eqnarray*}
\frac{S(0)}{N}\exp\Big((\eta N/2-1/2)T/N+\sup_{0\leq t\leq T/N}\widetilde{B}(t)\Big),\;N\in\nn 
\end{eqnarray*}
converges to the constant $m_\lambda e^{\eta T/2}$ in probability in the limit $N\rightarrow\infty$. Thus, the latter probability converges to $0$ in the limit $N\rightarrow\infty$ for all $\zeta$ small enough.

\bigskip

\noindent 4) To show that the second summand in expression (*) is zero, we first note that for every pair $0\leq s\leq t\leq T/N$ with $t-s\leq\zeta/N$ there is a $k\in\nn$ with $s,t\in[k\zeta/N,(k+2)\zeta/N]$. We use this observation and the union bound to conclude
\begin{eqnarray*}
&&\pp\Big(\sup_{0\leq s\leq t\leq T/N,t-s\leq\zeta/N}|M(t)-M(s)|>\Delta/2\Big)\\
&\leq&\sum_{k=0}^{\lfloor T/\zeta\rfloor-1}\pp\Big(\sup_{k\zeta/N\leq t\leq (k+2)\zeta/N} |M(t)-M(k\zeta/N)|>\Delta/4\Big),
\end{eqnarray*}
where $\lfloor.\rfloor$ denotes the function taking a real number to its integer part.

Next, we use \eqref{Sexplicitly} to compute 
\begin{equation}\label{Smoments}\begin{split}
&\ev[S(t)]=\ev[S(0)]\exp(\eta Nt/2),\\ 
&\ev[S(t)^2]=\ev[S(0)^2]\exp((\eta N+1)t).
\end{split}\end{equation}
The inequality $f'_r(X_i(t))^2 X_i(t)S(t)\leq \sup_{x\in\rr}|f'_r(x)|^2 S(t)^2$, Fubini's Theorem and \eqref{Smoments} imply that the process $M(t)$, $t\geq0$ is a martingale. Applying the $L^2$-version of Doob's maximal inequality for non-negative submartingales we obtain
\begin{eqnarray*}
&&\sum_{k=0}^{\lfloor T/\zeta\rfloor-1}\pp\Big(\sup_{k\zeta/N\leq t\leq (k+2)\zeta/N} |M(t)-M(k\zeta/N)|>\Delta/4\Big)\\
&\leq&\frac{16}{\Delta^2} \sum_{k=0}^{\lfloor T/\zeta\rfloor-1}\ev[(M((k+2)\zeta/N)-M(k\zeta/N))^2].
\end{eqnarray*}
By the Ito isometry the latter expression can be computed to 
\begin{eqnarray*}
&&\frac{16}{\Delta^2 N^2} \sum_{k=0}^{\lfloor T/\zeta\rfloor-1} \sum_{i=1}^N 
\ev\Big[\int_{k\zeta/N}^{(k+2)\zeta/N} f_r'(X_i(u))^2X_i(u)S(u)\;du\Big]\\
&\leq&\frac{32\sup_{x\in\rr}|f'_r(x)|^2}{\Delta^2N^2}\ev\Big[\int_0^{(\lfloor T/\zeta\rfloor+1)\zeta/N} S(u)^2\;du\Big].
\end{eqnarray*}
By Fubini's Theorem (note that the integrand is non-negative) and \eqref{Smoments} one deduces that the right-hand side is equal to 
\begin{eqnarray*}
\frac{32\sup_{x\in\rr}|f'_r(x)|^2\ev[S(0)^2]}{\Delta^2 N^2(\eta N+1)}\Big(\exp\Big((\eta N+1)(\lfloor T/\zeta\rfloor+1)\zeta/N\Big)-1\Big).
\end{eqnarray*}
Using \eqref{m_lambda} we conclude that the latter expression tends to $0$ in the limit $N\rightarrow\infty$ for any fixed $\zeta>0$. \ep

\subsection{Identification of the limit point}

In this section we will uniquely characterize the limit points of the sequence $Q^N_T$, $N\in\nn$ and thereby complete the proof of Theorem \ref{main_result} (a). To this end, we fix a convergent subsequence $Q^{N_k}_T$, $k\in\nn$ of the sequence $Q^N_T$, $N\in\nn$ and let $Q^\infty_T$ be its limit. By the Skorokhod Representation Theorem in the form of Theorem 3.5.1 in \cite{du}, there exist $C([0,T],M_1(\rr))$-valued random variables $\widetilde{\rho}_k$, $k\in\nn$ with laws $Q^{N_k}_T$, $k\in\nn$ converging to a limiting random variable $\widetilde{\rho}_\infty$ with law $Q^\infty_T$ almost surely. We start with the following lemma.

{\lemma\label{averageprocess} The sequence of functions $t\mapsto(\widetilde{\rho}_k(t),x)$, $k\in\nn$ converges in the space $C([0,T],\rr)$ to $t\mapsto m_\lambda e^{\eta t/2}$ in probability in the limit $k\rightarrow\infty$.\\\quad\\}
{\it Proof.} First, we fix a $k\in\nn$ and with a minor abuse of notation write $S(t)$ for $X_1(t)+\dots+X_{N_k}(t)$. By equation \eqref{Sexplicitly}, we have for every $\epsilon>0$:
\begin{eqnarray*}
&&\pp\Big(\sup_{0\leq t\leq T} \Big|(\widetilde{\rho}_k(t),x)-m_\lambda e^{\eta t/2}\Big|>\epsilon\Big)\\
&=&\pp\Big(\sup_{0\leq t\leq T} \Big|\frac{S(0)}{N_k}\exp\Big(\eta/2-1/(2N_k))t+N_k^{-1/2}W(t)\Big)-m_\lambda e^{\eta t/2}\Big|>\epsilon\Big)\\
&\leq&\pp\Big(\sup_{0\leq t\leq T} \Big|\frac{S(0)}{N_k}\exp\Big(-t/(2N_k)+N_k^{-1/2}W(t)\Big)-m_\lambda\Big|>\epsilon e^{-\eta T/2}\Big),
\end{eqnarray*} 
where $W(t)=N_k^{1/2}\widetilde{B}(t/N_k)$, $t\geq0$ is a standard Brownian motion. By Girsanov's Theorem the process $\frac{S(0)}{N_k}\exp\Big(-t/(2N_k)+N_k^{-1/2}W(t)\Big)-m_\lambda$, $t\geq0$ is a martingale. Hence, we can apply the $L^2$-version of Doob's maximal inequality for non-negative submartingales to estimate the latter upper bound from above by
\begin{eqnarray*}
&&\frac{e^{\eta T}}{\epsilon^2} \ev\Big[\Big(\frac{S(0)}{N_k}\exp\Big(-T/(2N_k)+N_k^{-1/2}W(T)\Big)-m_\lambda\Big)^2\Big]\\
&=&\frac{e^{\eta T}}{\epsilon^2}\Big(\frac{\ev[S(0)^2]}{N_k^2}\ev\Big[\exp\Big(-\frac{T}{N_k}+\frac{2}{\sqrt{N_k}}W(T)\Big)\Big]
-2m_\lambda\frac{\ev[S(0)]}{N_k}+m_\lambda^2\Big)\\
&=&\frac{e^{\eta T}}{\epsilon^2}\Big(\frac{\ev[S(0)^2]}{N_k^2}e^{T/N_k}-2m_\lambda\frac{\ev[S(0)]}{N_k}+m_\lambda^2\Big).
\end{eqnarray*}
It follows from \eqref{m_lambda} that the latter expression tends to $0$ in the limit $k\rightarrow\infty$. This finishes the proof of the lemma. \ep

\bigskip

Now, we are ready to prove that the system \eqref{pde}, \eqref{pdeic} must hold in the distributional sense almost surely under any limit point of the sequence $Q^N_T$, $N\in\nn$. 

{\prop\label{convergence} Let $Q^\infty_T$ be the limit of a convergent subsequence $Q^{N_k}_T$, $k\in\nn$ of the sequence $Q^N_T$, $N\in\nn$. Then under $Q^\infty_T$ the system \eqref{pde}, \eqref{pdeic} is satisfied in the distributional sense almost surely. Moreover, if $\widetilde{\rho}_\infty$ is a random variable with law $Q^\infty_T$, then it holds $\widetilde{\rho}_\infty\in C([0,T],M_1([0,\infty))$ with probability $1$.\\\quad\\}
{\it Proof.} 1) Fix a $Q^\infty_T$ as in the statement of the proposition. Let $\{g_1,g_2,\dots\}$ be a dense subset of the space ${\mathcal S}(\rr)$ of Schwartz functions on $\rr$ with respect to the topology of uniform convergence of functions and their first and second derivatives, such that each $g_r$ is infinitely differentiable with compact support. We claim that in order to prove the first assertion of the proposition it suffices to show that
\begin{eqnarray}
&&(\rho(t),g_r)-(\rho(0),g_r)=m_\lambda\int_0^t e^{\frac{\eta s}{2}}\Big(\rho(s),\frac{\eta}{2}g_r'+\frac{1}{2}xg_r''\Big)\;ds \label{inteq},\\
&&\rho(0)=\lambda \label{inteqic}
\end{eqnarray}
holds for all $r\in\nn$ and $t\in[0,T]$ with probability $1$ under $Q^\infty_T$. Indeed, this would imply that the system \eqref{inteq}, \eqref{inteqic} is satisfied for all $g\in{\mathcal S}(\rr)$ and all $t\in[0,T]$ with probability $1$ under $Q^\infty_T$. This would yield the first assertion of the proposition.  

\bigskip

\noindent 2) Since a countable union of null sets is a null set, it is enough to show that the system \eqref{inteq}, \eqref{inteqic} is satisfied for a fixed function $g_r$ and all $t\in[0,T]$ with probability $1$ under $Q^\infty_T$. It is clear that equation \eqref{inteqic} is satisfied with probability $1$ due to Assumption \ref{ass}. In order to show that equation \eqref{inteq} holds, we use Ito's formula to compute
\begin{eqnarray*}
d\frac{1}{N}\sum_{i=1}^N g_r(X_i(t))=\frac{1}{N}\sum_{i=1}^N \Big(\frac{\eta}{2}g_r'(X_i(t))S(t)+\frac{1}{2}g_r''(X_i(t))X_i(t)S(t)\Big)dt\\
+\frac{1}{N}\sum_{i=1}^N g_r'(X_i(t))\sqrt{X_i(t)S(t)} dW_i(t).
\end{eqnarray*}
Thus, one has the dynamics
\begin{eqnarray}
\;\;\;\;\;\;\;\;\;\;\;\;d(\rho^N(t),g_r)=D_r(t)\;dt+\frac{1}{N^{3/2}}\sum_{i=1}^N g_r'(Y_i(t))\sqrt{Y_i(t)S^Y(t)} dB_i(t),
\end{eqnarray}
where 
\begin{eqnarray*}
D_r(t)&=&\frac{1}{N^2}\sum_{i=1}^N \Big(\frac{\eta}{2}g_r'(Y_i(t))S^Y(t)+\frac{1}{2}g_r''(Y_i(t))Y_i(t)S^Y(t)\Big)\\
&=&\Big(\rho^N(t),\frac{\eta}{2}g_r'+\frac{1}{2}xg_r''\Big)\cdot(\rho^N(t),x),
\end{eqnarray*}
$t\geq0$, and $B_i(t)=N^{1/2} W_i(t/N)$, $1\leq i\leq N$, $S^Y(t)=Y_1(t)+\dots+Y_N(t)$ as before. Moreover, the inequality
\begin{eqnarray}
|g_r'(Y_i(t))|\sqrt{Y_i(t)S^Y(t)}\leq\sup_{x\in\rr}|g_r'(x)| S^Y(t),
\end{eqnarray}
Fubini's Theorem and \eqref{Smoments} show that the process
\begin{eqnarray}
\;\;\;\;\;\;M_r(t)=\frac{1}{N^{3/2}}\sum_{i=1}^N \int_0^t g_r'(Y_i(u))\sqrt{Y_i(u)S^Y(u)} dB_i(u),\;t\geq0
\end{eqnarray}
is a martingale. Hence, using the $L^2$-version of Doob's maximal inequality for non-negative submartingales and the Ito isometry, we obtain for every $\epsilon>0$:
\begin{eqnarray*}
&&\pp\Big(\sup_{0\leq t\leq T} \Big|(\rho^N(t),g_r)-(\rho^N(0),g_r)-\int_0^t D_r(u)\;du\Big|>\epsilon\Big)\\
&&\leq\epsilon^{-2}\ev[M_r(T)^2]\\
&&=\frac{1}{\epsilon^2 N^3}\sum_{i=1}^N \ev\Big[\int_0^T g_r'(Y_i(u))^2 Y_i(u)S^Y(u) du\Big]\\
&&\leq\frac{\sup_{x\in\rr}|g_r'(x)|^2}{\epsilon^2 N^3}\ev\Big[\int_0^T S^Y(u)^2\;du\Big]. 
\end{eqnarray*} 
Using $S^Y(u)=S(u/N)$, $u\geq0$, Fubini's Theorem and \eqref{Smoments} we can compute the latter upper bound to
\begin{eqnarray}
\frac{\sup_{x\in\rr}|g_r'(x)|^2}{\epsilon^2 N^3}\cdot\frac{N\;\ev[S(0)^2]}{\eta N+1}\Big(e^{(\eta N+1)T/N}-1\Big).
\end{eqnarray}
This expression tends to $0$ in the limit $N\rightarrow\infty$ for every $\epsilon>0$ due to \eqref{m_lambda}.

\bigskip

\noindent 3) Next, we recall the definition of the random variables $\widetilde{\rho}_k$, $k\in\nn$ and $\widetilde{\rho}_\infty$ prior to Lemma \ref{averageprocess}. In view of the latter, we may and will assume that the sequence of functions $t\mapsto(\widetilde{\rho}_k(t),x)$, $k\in\nn$ converges to $t\mapsto m_\lambda e^{\eta t/2}$ in the space $C([0,T],\rr)$ with probability $1$ (otherwise we pass to a suitable subsequence). It follows that the random variables $\Xi_k$ given by
\begin{eqnarray*}
\sup_{0\leq t\leq T} \Big|(\widetilde{\rho}_k(t),g_r)-(\widetilde{\rho}_k(0),g_r)
-\int_0^t (\widetilde{\rho}_k(u),\eta/2 g_r'+x/2g_r'')\cdot(\widetilde{\rho}_k(u),x)\;du\Big|,
\end{eqnarray*}
converge almost surely in the limit $k\rightarrow\infty$ to
\begin{eqnarray*}
\sup_{0\leq t\leq T} \Big|(\widetilde{\rho}_\infty(t),g_r)-(\widetilde{\rho}_\infty(0),g_r)
-\int_0^t (\widetilde{\rho}_\infty(u),\eta/2 g_r'+x/2g_r'')\cdot m_\lambda e^{\eta u/2}\;du\Big|,
\end{eqnarray*}
which we call $\Xi_\infty$. Finally, using the Portmanteau Theorem and the final result of step 2 we obtain for every $\epsilon>0$:
\begin{eqnarray*}
&&\pp(\Xi_\infty>\epsilon)\leq\liminf_{k\rightarrow\infty}\pp(\Xi_k>\epsilon)\\
&&=\liminf_{k\rightarrow\infty}\pp\Big(\sup_{0\leq t\leq T} \Big|(\rho^{N_k}(t),g_r)-(\rho^{N_k}(0),g_r)-\int_0^t D_r(u)\;du\Big|>\epsilon\Big)=0.
\end{eqnarray*}
Since the law of $\widetilde{\rho}_\infty$ is given by $Q^\infty_T$, it follows that equation \eqref{inteq} holds $Q^\infty_T$-almost surely. 

\bigskip

\noindent 4) To prove the second assertion of the proposition, we note that $\widetilde{\rho}_k(t)([0,\infty))=1$ holds for all $k\in\nn$ and $t\in[0,T]$ almost surely. This is a consequence of the representation of $X_1,\dots,X_N$ as time-changed squared Bessel processes (see equations (12.7)-(12.9) in \cite{fk}) and the properties of the latter (see e.g. chapter XI of \cite{ry}). Thus, the Portmanteau Theorem implies $\widetilde{\rho}_\infty(t)([0,\infty))=1$ for all $t\in[0,T]$ on the same set of full probability. \ep 

\bigskip

In view of Propositions \ref{tightness} and \ref{convergence}, the proof of Theorem \ref{main_result} (a) is complete, once we show that the solution of the Cauchy problem 
\begin{eqnarray}
\;\;\;\;\;\;&&\forall g\in{\mathcal S}(\rr):\;
(\rho(t),g)-(\rho(0),g)=m_\lambda\int_0^t e^{\frac{\eta s}{2}}\Big(\rho(s),\frac{\eta}{2}g'+\frac{1}{2}xg''\Big)\;ds \label{inteq2},\\
\;\;\;\;\;\;&&\rho(0)=\lambda \label{inteq2ic}
\end{eqnarray}
in $C([0,T],M_1([0,\infty)))$ is unique. 

{\prop\label{uniqueness} The solution of the Cauchy problem \eqref{inteq2}, \eqref{inteq2ic} in the space $C([0,T],M_1([0,\infty))$ is unique.\\\quad\\}
{\it Proof.} 1) Let $\mu,\nu\in C([0,T],M_1([0,\infty)))$ be two solutions of the problem \eqref{inteq2}, \eqref{inteq2ic}. Moreover, define the operator 
\begin{eqnarray}
Lh=\frac{\partial h}{\partial t}-m_\lambda e^{\frac{\eta t}{2}}\frac{\eta}{2}\frac{\partial h}{\partial x}
-m_\lambda e^{\frac{\eta t}{2}}\frac{x}{2}\frac{\partial^2 h}{\partial x^2}
\end{eqnarray}
acting on the space $C^{1,2}_c([0,T]\times\rr)$ of continuous real-valued functions on $[0,T]\times\rr$ having compact support in $(0,T)\times\rr$, a continuous time derivative and two continuous spatial derivatives. For the same reason as in the remark preceeding Theorem A.1 in \cite{ga}, the problem \eqref{inteq2}, \eqref{inteq2ic} is equivalent to the problem  
\begin{equation}\label{inteq3}\begin{split}
\forall h\in C^{1,2}_c([0,T]\times\rr):(\rho(t),h(t,.))-(\rho(0),h(0,.))=\int_0^t (\rho(s),(Lh)(s,.))ds
\end{split}\end{equation}
with initial condition $\rho(0)=\lambda$. This follows by an approximation of functions $h\in C^{1,2}_c([0,T]\times\rr)$ by functions on $[0,T]\times\rr$ which are Schwartz functions in $x$ for every fixed $t$ and piecewise constant in $t$. Since for all $h\in C^{1,2}_c([0,T]\times\rr)$ it holds $(\mu(T),h(T,.))=(\nu(T),h(T,.))=0$, we conclude from \eqref{inteq3}: 
\begin{equation}\begin{split}
\forall h\in C^{1,2}_c([0,T]\times\rr):\;\int_0^T (\mu(s),(Lh)(s,.))\;ds=\int_0^T (\nu(s),(Lh)(s,.))\;ds.
\end{split}\end{equation}

\noindent 2) In view of the latter equation, it suffices to show that the space 
\begin{eqnarray*}
Lh,\;h\in C^{1,2}_c([0,T]\times\rr) 
\end{eqnarray*}
contains all functions on $[0,T]\times\rr$ of the form $(t,x)\mapsto l(t)r(x)$, where $l$ is an infinitely differentiable function on $[0,T]$ with compact support contained in $[t_0,t_1]$ for some $0<t_0<t_1<T$ and $r$ is an infinitely differentiable function on $\rr$ with compact support contained in $[a,b]$ for some $0<a<b$. To this end, we first let ${\mathcal R}$ be a bounded closed rectangle of the form $[\widetilde{t}_0,\widetilde{t}_1]\times[\widetilde{a},\widetilde{b}]$ which is contained in $(0,T)\times\rr$ and such that $[t_0,t_1]\times[a,b]$ is contained in the interior of ${\mathcal R}$. By Theorem 9 in section 1.5 of \cite{fr} there exists a function $h:{\mathcal R}\rightarrow\rr$ such that the derivatives $\frac{\partial h}{\partial t}$, $\frac{\partial h}{\partial x}$, $\frac{\partial^2 h}{\partial x^2}$ exist and are continuous on ${\mathcal R}$ and it holds $(Lh)(t,x)=l(t)r(x)$ in the strong sense on ${\mathcal R}$. 

Next, we apply the maximum principle in the form of Lemma 4 in section 2.1 of \cite{fr} on the rectangles
\begin{eqnarray*}
[\widetilde{t}_0,\widetilde{t}_1]\times[\widetilde{a},a],\;[t_1,\widetilde{t}_1]\times[\widetilde{a},\widetilde{b}],
\;[\widetilde{t}_0,\widetilde{t}_1]\times[b,\widetilde{b}],\;[\widetilde{t}_0,t_0]\times[\widetilde{a},\widetilde{b}] 
\end{eqnarray*}
to the functions $h$ and $-h$ to conclude that the function $h$ is constant on the complement of $[t_0,t_1]\times[a,b]$ in ${\mathcal R}$. Hence, without loss of generality we may assume that $h$ is constantly equal to zero on the complement of $[t_0,t_1]\times[a,b]$ in ${\mathcal R}$ (otherwise we subtract a constant from $h$). Thus, we can extend $h$ to a function $\widetilde{h}\in C^{1,2}_c([0,T]\times\rr)$ such that $\widetilde{h}=h$ on ${\mathcal R}$ and $\widetilde{h}=0$ on the complement of ${\mathcal R}$. It follows that $(L\widetilde{h})(t,x)=l(t)r(x)$ holds in the strong sense on $[0,T]\times\rr$, which finishes the proof. \ep  

\subsection{Stochastic representation of the limit point}

In this subsection we prove part (b) of Theorem \ref{main_result}.\\\quad\\
{\it Proof of Theorem \ref{main_result} (b).} In view of Proposition \ref{uniqueness} it suffices to show that the one-dimensional distributions $\xi(t)$, $t\in[0,T]$ of the process $\Gamma(t)=Z\Big(\int_0^t e^{\eta s/2}m_\lambda\;ds\Big)$, $t\in[0,T]$ form a distributional solution to the Cauchy problem \eqref{pde}, \eqref{pdeic} in $C([0,T],M_1([0,\infty)))$. It is clear that $\xi(t)$, $t\in[0,T]$ satisfies \eqref{pdeic} and that it is an element of $C([0,T],M_1([0,\infty)))$, since the squared Bessel process $Z$ has continuous paths and takes values in $[0,\infty)$ (see e.g. chapter XI in \cite{ry}). To prove that $\xi(t)$, $t\in[0,T]$ satisfies \eqref{pde}, we use the time-change formalism for Brownian motion (see e.g. \cite{ok}, chapter 8.5) to deduce
\begin{eqnarray*}
d\Gamma(t)&=&\frac{d}{dt}\Big(\int_0^t e^{\eta s/2}m_\lambda\;ds\Big)\frac{\eta}{2}\;dt
+\sqrt{\frac{d}{dt}\Big(\int_0^t e^{\eta s/2}m_\lambda\;ds\Big)}\sqrt{\Gamma(t)}\;d\widetilde{\beta}(t)\\
&=&e^{\eta t/2}m_\lambda\frac{\eta}{2}\;dt+e^{\eta t/4}\sqrt{m_\lambda}\sqrt{\Gamma(t)}\;d\widetilde{\beta}(t),
\end{eqnarray*}
$t\geq0$, where $\widetilde{\beta}$ is an appropriate standard Brownian motion. Hence, Ito's formula shows 
\begin{eqnarray*}
g(\Gamma(t))-g(\Gamma(0))&=&\int_0^t m_\lambda e^{\eta s/2}\Big(\frac{\eta}{2}g'(\Gamma(s))+\frac{1}{2}\Gamma(s)g''(\Gamma(s))\Big)\;ds\\
&+&\int_0^t g'(\Gamma(s))e^{\eta s/4}\sqrt{m_\lambda}\sqrt{\Gamma(s)}\;d\widetilde{\beta}(s)
\end{eqnarray*}
for all $g\in{\mathcal S}(\rr)$ and $t\in[0,T]$. Moreover, the expectations $\ev[\Gamma(t)]$, $t\in[0,T]$ are finite and uniformly bounded (this is evident from the definition of squared Bessel processes in chapter XI of \cite{ry} if $2\eta$ is an integer; in the general case, this can be deduced by comparing $Z$ with a squared Bessel process of integer index $2\eta'>2\eta$ using, for example, the comparison theorems of section 3 in chapter IX of \cite{ry}). Thus, by taking the expectation in the latter equation and applying Fubini's Theorem we get
\begin{eqnarray}
\;\;\;\;\;\;(\xi(t),g)-(\xi(0),g)=\int_0^t \Big(\xi(s),m_\lambda e^{\eta s/2}\Big(\frac{\eta}{2}g'+\frac{1}{2}xg''\Big)\Big)\;ds
\end{eqnarray} 
for all $g\in{\mathcal S}(\rr)$ and $t\in[0,T]$. Thus, $\xi(t)$, $t\in[0,T]$ solves the equation \eqref{pde} in the distributional sense. \ep 

\section*{Acknowledgements}
The author thanks Amir Dembo for his comments on this work. He is also grateful to Adrian Banner, Robert Ferholz, Ioannis Karatzas and Vassilios Papathanakos who brought the question on the large $N$ behavior of volatility-stabilized market models to his attention.

\end{document}